\documentclass[11pt]{article}
\usepackage[a4paper,margin=28mm]{geometry}
\usepackage{amsmath,amssymb,amsthm,mathtools}
\usepackage{mathrsfs}
\usepackage[hidelinks]{hyperref}
\usepackage[shortlabels]{enumitem}
\usepackage{dsfont}
\numberwithin{equation}{section}
\newtheorem{theorem}{Theorem}[section]
\newtheorem{proposition}[theorem]{Proposition}
\newtheorem{lemma}[theorem]{Lemma}
\theoremstyle{remark}

\newcommand{\N}{\mathbb N_0}
\newcommand{\C}{\mathbb C}
\newcommand{\R}{\mathbb R}

\newcommand{\dd}{\mathrm d}
\newcommand{\e}{\mathrm{e}}
\usepackage{upgreek}
\newcommand{\iu}{\mathrm{i}}
\newcommand{\mfrak}{\mathfrak}
\newcommand{\scr}{\mathscr}
\renewcommand{\cal}{\mathcal}

\newcommand{\norm}[1]{\left\lVert#1\right\rVert}
\newcommand{\ind}[1]{{\mathds{1}_{{#1}}}}
\mathtoolsset{showonlyrefs}

\DeclareMathOperator{\supp}{supp}
\renewcommand{\Re}{\operatorname{Re}}

\title{Group-free radial Kunze--Stein phenomenon\protect\\on affine buildings}
\author{Łukasz Pawlak}
\date{}
\newcommand{\ground}{\boldsymbol{\Phi}}
\begin{document}
\maketitle
\begin{abstract}
    We give a geometric proof of the radial Kunze--Stein phenomenon on the
    good vertices of a thick, locally finite, regular affine building of
    irreducible type. The proof uses a positive spherical function, its
    joint eigenvalues for sphere averages, and the growth of vector-distance
    spheres. Conjugation by this function and complex interpolation yield
    a weighted $\ell^2$ estimate and the boundedness of radial operators with
    $\ell^p$ kernels for $1\leq p<2$. We also present proofs by the
    averaging operators norms, by the Schur test
    and by Helgason transform. A further interpolation gives weighted
    mixed-exponent estimates. No transitivity assumption on the automorphism
    group is required; in particular the argument applies to buildings of
    type $\widetilde A_2$ without a Bruhat--Tits realization.
\end{abstract}
\let\thefootnote\relax\footnote{Project with registration number 2023/49/B/ST1/00678/R was
financed by the National Science Centre (NCN).}

\section{Introduction}\label{sec:intro}
Young's inequality bounds the $L^2$ norm of convolution by the $L^1$
norm of its kernel. On certain noncompact groups the kernel may instead
belong to $L^p$, for any $1\leq p<2$. This improvement is called the
Kunze--Stein phenomenon. It was established for $\mathrm{SL}(2,\R)$ in
\cite{KunzeStein}. Herz's majorization principle
\cite{Herz1970} relates such estimates to a distinguished positive
spherical function. Cowling established the Kunze--Stein phenomenon for
connected real semisimple Lie groups with finite center
\cite{Cowling1978}; see also his exposition \cite{CowlingHerz1997}.

These estimates also occur in the study of functions of Laplace operators.
For noncompact symmetric spaces, Cowling, Giulini and Meda obtain
$L^p$--$L^q$ mapping estimates for the heat semigroup in
\cite[Theorem 3.2]{CowlingGiuliniMeda1993}. The corresponding
continuous-time heat-semigroup estimates on homogeneous trees are given
in \cite[Theorem 2.2]{CowlingMedaSetti2000}. The Kunze--Stein phenomenon
on the full isometry group of a homogeneous tree is treated in
\cite{VecaTree}. For applications to dispersive equations and the
weighted mixed-exponent interpolation that motivates our last result,
see \cite[Lemma 5.1]{AnkerPierfeliceVallarino}.

Affine buildings provide a setting in which these estimates can be
studied independently of a transitive group action. Their apartments
are Euclidean Coxeter complexes, and vector distance replaces ordinary
distance on a tree. Macdonald developed spherical harmonic analysis
for groups of $p$-adic type \cite{Macdonald1971}. For buildings of type
$\widetilde A_r$, see also Cartwright \cite{Cartwright2001}.
The geometric theory for arbitrary locally finite regular affine
buildings, including those not arising from groups, is developed in
\cite{ParkinsonHecke,ParkinsonSHA}. We use the notation of the recent
works of R\'emy--Trojan and Papageorgiou--Trojan
\cite{RemyTrojan,PapageorgiouTrojan}. In the Bruhat--Tits setting, the
radial Herz and Kunze--Stein inequalities are recalled in
\cite[Appendix A]{PapageorgiouTrojan}.

Here are our main results. For a fixed good vertex $o$, a radial function $k$
is determined by its profile $\tilde k$ through
$k(x)=\tilde k(\sigma(o,x))$, where $\sigma$ denotes the vector distance.
We use the following convolution convention:
\begin{equation}
    \label{eq:conv}
    T_ku(x)=(u\times k)(x)
    =\sum_{y\in V_g}u(y)\tilde k(\sigma(x,y)).
\end{equation}
The sum is initially defined for finitely supported $u$.
Let $P_\lambda(z)$ be the Macdonald spherical function defined
in section \ref{sec:prelim}, and $\varphi_0(x) = P_{\sigma(o, x)}(0)$
be the ground spherical function.
\begin{theorem}[Radial Kunze--Stein estimate]\label{thm:intro-KS}
    Let $\mathscr X$ be a thick, locally finite, regular irreducible affine
    building.
    For a radial complex function $k(x)=\tilde k(\sigma(o,x))$, put
    \[
        \mathcal H(|k|)(0):=\sum_{x\in V_g}|k(x)|\varphi_0(x)
        =\sum_{\lambda\in P^+}N_\lambda|\tilde k(\lambda)|P_\lambda(0).
    \]
    If this quantity is finite, the operator
    $T_ku(x)$ defined in \eqref{eq:conv} initially on finite-support
    functions, extends boundedly to $\ell^2(V_g)$, and
    \[
        \|T_ku\|_2\leq\mathcal H(|k|)(0)\|u\|_2, \qquad
        u \in \ell^2(V_g).
    \]
    In particular, for $1\leq p<2$ and radial $k\in\ell^p(V_g)$,
    \[
        \|T_ku\|_2\leq\|\varphi_0\|_{p'}\|k\|_p\|u\|_2,
        \qquad u\in\ell^2(V_g),\quad \frac 1 p + \frac 1 {p'} = 1.
    \]
\end{theorem}
\noindent
The notation $\mathcal H(|k|)(0)$ is the Helgason transform of the
radial function $|k|$ at zero; its geometric definition is recalled
in the preliminaries.

We then interpolate the result of Theorem \ref{thm:intro-KS};
for $2\leq a\leq b<\infty$, put $1/s=1/a+1/b$ and $1/a+1/a'=1$.
If $k$ is radial and $\varphi_0^{2/b}k\in\ell^s(V_g)$, then $T_k$
extends to a bounded map $\ell^{a'}(V_g)\to\ell^b(V_g)$ satisfying
\[
    \|T_ku\|_b\leq\|u\|_{a'}\|\varphi_0^{2/b}k\|_s, \qquad
    u \in \ell^2(V_g).
\]
The weighted $\ell^2$ inequality already follows from the known
sphere-operator norm formula
$\|A_\lambda\|_{2\to2}=P_\lambda(0)$
\cite[Theorem 6.3]{ParkinsonSHA}. Our purpose is to exhibit its geometric
mechanism and the resulting interpolation estimates.

The interpolated Kunze--Stein estimate allows one to deduce
$\ell^{a'}$--$\ell^b$ bounds for radial convolution operators directly
from pointwise estimates of their kernels. Indeed, if
$|k(x)|\leq M(x)$, then
\[
    \|T_k\|_{a'\to b} \leq
    \left( \sum_{x \in V_g}
    M(x)^s \varphi_0(x)^{2s/b} \right)^{1/s},
    \qquad 2\leq a\leq b<\infty, \quad \frac 1 s = \frac 1 a+\frac 1 b.
\]
In particular, pointwise bounds for radial heat or Schrödinger
kernels yield corresponding smoothing or dispersive estimates
whenever the sum is finite, with its dependence on time
determining the resulting operator-norm bound.

Section~\ref{sec:prelim} collects the geometric, spherical, spectral and
interpolation facts used in the proofs. Section~\ref{sec:kunze} gives
four proofs of the radial estimate: averaging operator norm one,
conjugation and interpolation,
the weighted Schur test, and the Helgason transform argument communicated to us
by Effie Papageorgiou. After that, it contains the mixed-exponent
estimate proof.

\section{Preliminaries}\label{sec:prelim}
Throughout, let $\mathscr X$ be a thick, locally finite, regular affine
building of irreducible type. Thus each panel of type $i$ lies in
$q_i+1$ chambers, with $1<q_i<\infty$ independent of the panel.
An apartment is a copy of a Euclidean Coxeter complex; any two
simplices lie in an apartment, and apartment identifications can be
chosen to fix their intersection. A panel is a face of codimension one;
its type specifies the corresponding simple affine reflection.

Let $\Phi\subset\mathfrak a$ be the associated finite root system,
possibly nonreduced, where $\mathfrak a$ is an $r$-dimensional Euclidean
space. Choose positive roots $\Phi^+$ with simple roots
$\alpha_1,\ldots,\alpha_r$, and set
$\Phi^{++}=\{\alpha\in\Phi^+:\alpha/2\notin\Phi\}$.
For a positive root $\alpha \in \Phi^+$, its height is defined as
\begin{equation}
    \operatorname{ht}(\alpha) = \sum_{i = 1}^r n_i
\end{equation}
whenever $\alpha = \sum_{i = 1}^r n_i \alpha_i$ with $n_i \in \N$.
As $\Phi$ is irreducible, there exists a unique highest root
$\alpha_0 = \sum_{i = 1}^r m_i \alpha_i$ with $m_i \in \N$.
For $\alpha\in\Phi$, write
$\alpha^\vee=2\alpha/\langle\alpha,\alpha\rangle$.
The coroot and coweight lattices and the dominant cone are
\[
    Q=\sum_{\alpha\in\Phi}\mathbb Z\alpha^\vee,\qquad
    P=\bigoplus_{i=1}^r\mathbb Z\lambda_i,\qquad
    P^+=\bigoplus_{i=1}^r\mathbb N_0\lambda_i,
\]
where $\lambda_i$'s form basis dual to $\alpha_j$'s:
$\langle\lambda_i,\alpha_j\rangle=\delta_{ij}$.
Lattice dual to $Q$ is labeled $L$ and called the weight lattice.
Denote the simple
affine reflections by $r_0,\ldots,r_r$, that is set
$r_0 = r_{\alpha_0; 1}$ and $r_i = r_{\alpha_i; 0}$, where
\begin{equation}
    r_{\alpha; k}(x) = x - \left( \langle x, \alpha \rangle - k\right) \alpha^\lor.
\end{equation}
The finite Weyl group is one generated by $r_i, i = 1, \ldots, r$ and denoted $W$.
Its longest element is $w_0$, and
$\lambda^*=-w_0\lambda$ for $\lambda\in P^+$.
The affine and extended affine Weyl groups are
$W^a=Q\rtimes W$ and $\widetilde W^a=P\rtimes W$, with the natural
reflection action on the translation lattices. 
The fundamental sector in $\mfrak a$ defined as
\begin{equation}
    S_0 = \{x \in \mfrak a \colon
    \langle x, \alpha_i \rangle \geq 0, i = 1, \ldots, r\}
    = \bigoplus_{i = 1}^r \R_{\geq 0}\lambda_i
\end{equation}
is the fundamental domain for the action of $W$ on $\mfrak a$.
The fundamental chamber (or fundamental alcove) is defined by
\begin{equation}
    C_0 = \{x \in \mfrak a \colon \langle x, \alpha_0 \rangle \leq 1,
    \langle x, \alpha_i \rangle \geq 0,  i = 1, \ldots, r \}.
\end{equation}
Let $\{\lambda + wS_0\colon \lambda \in P, w \in W^a\}$ be the
set of sectors in $\mfrak a$. $\lambda$ is called the base vertex of the sector.
Subset of an apartment of the building
$\scr X$ is called a sector, when it belongs to just defined set of sectors
in $\mfrak a$ under identification of the apartment with model Coxeter complex.

The good vertices $V_g$ are those represented by $P$ in apartment
identification with the model apartment.
For reduced $\Phi$ they are all special vertices. In type
$\widetilde C_r$, use $\Phi=C_r$ if $q_0=q_r$, and $\Phi=BC_r$
otherwise. In the latter case only the type-$0$ special vertices are
good. The same convention includes semiregular trees via $BC_1$;
see \cite[Sections 1.4 and 1.7]{ParkinsonSHA}.
For $x,y\in V_g$, their vector distance $\sigma(x,y)\in P^+$ is the
dominant representative of the displacement from $x$ to $y$ in the
orbit of group action by $W$. Set
\begin{equation}\label{eq:spheres}
    V_\lambda(x)=\{y\in V_g:\sigma(x,y)=\lambda\},\qquad
    N_\lambda=|V_\lambda(x)|, \qquad x \in V_g, \lambda \in P^+.
\end{equation}
As denoted, cardinalities $N_\lambda$ are independent of $x$; they are
finite as well. Define averaging operators
\begin{equation}\label{eq:averages}
    A_\lambda u(x)=\frac 1 {N_\lambda}\sum_{y\in V_\lambda(x)}u(y).
\end{equation}
We have
\begin{equation}\label{eq:opposition}
    \sigma(y,x)=\sigma(x,y)^*,\qquad
    N_{\lambda^*}=N_\lambda,\qquad A_\lambda^*=A_{\lambda^*}.
\end{equation}
Fix $o\in V_g$. For a radial function $k(x)=\tilde k(\sigma(o,x))$,
our convolution convention \eqref{eq:conv} gives
\begin{equation}\label{eq:kernel}
    T_ku(x)
    =\sum_{\lambda\in P^+}N_\lambda\tilde k(\lambda)A_\lambda u(x).
\end{equation}
The adjoint radial kernel has profile
$\widetilde{k^\dagger}(\lambda)=\overline{\tilde k(\lambda^*)}$.
For a profile $g$ define its radial lift $k_g(x)=g(\sigma(o,x))$ and
notice that
\[
    \|k_g\|_p=\left(\sum_{\lambda\in P^+}N_\lambda|g(\lambda)|^p\right)^{1/p},
    \qquad 1\leq p<\infty,
\]
(use the supremum convention for $p=\infty$).
We use $p'$ for the conjugate exponent, including $1'=\infty$.
We use counting measure on $V_g$,
sums over $\lambda$ run over $P^+$ and sums over $x$
and $y$ run over $V_g$ unless specified otherwise.

For $w=r_{i_1}\cdots r_{i_l}\in W^a$ reduced, put
$q_w=q_{i_1}\cdots q_{i_l}$. This is independent of the expression.
Extend this notion to $\widetilde W^a$ by assigning parameter $1$ to the stabilizer
of the fundamental alcove (the length-zero elements). For a finite
subset $D$ of $W$ write $D(q^{-1})=\sum_{w\in D}q_w^{-1}$, and let
$W_\lambda=\{w\in W:w\lambda=\lambda\}$.

\subsection{The spherical functions}
For an indivisible root $\alpha$, let $q_\alpha=q_i$ if
$\alpha\in W\alpha_i$, $1\leq i\leq r$. Put
\[
    R_1 = \{\alpha \in \Phi \colon 2\alpha \not\in \Phi\}, \qquad
    R_2 = \{\alpha \in \Phi \colon \alpha / 2 \not\in \Phi\}, \qquad
    R_3 = R_1 \cap R_2.
\]
Define
\[
    \tau_\alpha=\begin{cases}
        q_\alpha,&\alpha\in R_3, \\
        q_0,&\alpha \in R_1 \setminus R_3, \\
        q_\alpha/q_0,&\alpha \in R_2 \setminus R_3, \\
        1,&\alpha\notin\Phi.
    \end{cases}
\]
Put
\[
    \chi_0(v) 
    =\prod_{\alpha\in\Phi^+}\tau_\alpha^{\langle v,\alpha\rangle},
    \qquad v\in\mathfrak a.
\]
By \cite[Proposition 1.5 and Proposition A.1]{ParkinsonSHA}, 
\begin{equation}\label{eq:sphere}
    N_\lambda=\frac{W(q^{-1})}{W_\lambda(q^{-1})}\chi_0(\lambda),
    \qquad \chi_0(\lambda)=q_{t_\lambda},
    \qquad\lambda\in P^+.
\end{equation}
Here $t_\lambda$ is translation by $\lambda$.

We use additive spectral parameters $z\in\mathfrak a_\C$ and the
complex bilinear extension of the Euclidean inner product (or equivalently,
sesquilinear inner product linear in the first variable, since the second
variable will always come from $\mathfrak a$). Macdonald's
formula, in the building form of \cite[Section 2]{ParkinsonSHA}, is
\begin{align}
    P_\lambda(z)&=\frac{\chi_0(\lambda)^{-1/2}}{W(q^{-1})}
    \sum_{w\in W}\mathbf c(wz)\e^{\langle wz,\lambda\rangle},
    \label{eq:macdonald}
\end{align}
where the $\mathbf c$-function takes the form
\begin{align}
    \mathbf c(z)&=\prod_{\alpha\in\Phi^+}
    \frac{1-\tau_\alpha^{-1}\tau_{\alpha/2}^{-1/2}
    \e^{-\langle z,\alpha^\vee\rangle}}
    {1-\tau_{\alpha/2}^{-1/2}\e^{-\langle z,\alpha^\vee\rangle}}.
    \label{eq:c}
\end{align}
The group-theoretic formula is due to Macdonald
\cite[Theorem 4.1.2]{Macdonald1971}; the generalization to arbitrary regular
affine buildings is supplied by \cite{ParkinsonHecke,ParkinsonSHA}.
The whole sum extends to an exponential polynomial, the values of $P_\lambda(z)$
where the denominator of the $\mathbf c$-function equals zero are obtained
by taking proper limits. For example, for reduced $\Phi$,
\[
    \mathbf c(z)=\prod_{\alpha\in\Phi^+}
    \frac{1-q_\alpha^{-1}\e^{-\langle z,\alpha^\vee\rangle}}
    {1-\e^{-\langle z,\alpha^\vee\rangle}}.
\]
Define the spherical function by
\begin{equation}\label{eq:phi}
    \varphi_z(x)=P_{\sigma(o,x)}(z), \qquad z \in \mfrak a_\C, x \in V_g,
\end{equation}
and set
\begin{equation}
    \ground(\lambda)=P_\lambda(0), \qquad \lambda \in P^+.
\end{equation}
$\varphi_0(x)$ is called ground spherical function.
By \cite[Theorem 6.3]{ParkinsonSHA}, we have
\begin{equation}\label{eq:5}
    \norm{A_\lambda}_{2 \to 2} = \Phi(\lambda), \qquad \lambda \in P^+.
\end{equation}

\subsection{The geometric boundary formula}
The space of sectors is equipped with equivalence relation: sectors are
equivalent, when their intersection contains a sector.
Let $\Omega$ be the boundary of $\scr X$ --- the set of sector classes.
For $\omega \in \Omega$ and $x\in V_g$ write $[x, \omega]$ for
the sector with a base vertex $x$ in the equivalence class of $\omega$.
For $x \in V_g$ define probability harmonic measure
$\nu_x$, by assigning mass $N_\lambda^{-1}$
to each shadow $\Omega(x,y)=\{\omega:y\in[x,\omega]\}$ for
$y\in V_\lambda(x)$. 
The vector-valued Busemann function is
\[
    h(x,y;\omega)=\sigma(x,v)-\sigma(y,v),
    \qquad x, y \in V_g, \omega \in \Omega,
\]
where $v$ is any good vertex lying in common subsector of $[x,\omega]$
and $[y,\omega]$. It satisfies {cocycle identity}
$h(x,z;\omega)=h(x,y;\omega)+h(y,z;\omega)$.
For any $x, y \in V_g$ the measures $\nu_x, \nu_y$ are equivalent, with
Radon--Nikodym derivative
\begin{align}
    \frac{\dd\nu_y}{\dd\nu_x}(\omega) &=\chi_0(h(x,y;\omega));\label{eq:RN}
 \end{align}
see \cite[Theorem 3.17(ii)]{ParkinsonSHA}. Using harmonic measures one can
express Macdonald spherical function as
 \begin{align}
    P_\lambda(z)&=\int_\Omega
    \chi_0(h(x,y;\omega))^{1/2}\e^{\langle z,h(x,y;\omega)\rangle}
    \,\dd\nu_x(\omega),\quad y\in V_\lambda(x).
    \label{eq:boundary}
\end{align}
See \cite[(6.6)]{ParkinsonSHA}.
In particular,
\begin{equation}\label{eq:affinity}
    \varphi_0(x)=\int_\Omega\chi_0(h(o,x;\omega))^{1/2}\,\dd\nu_o(\omega)
    =\int_\Omega\left(\frac{\dd\nu_x}{\dd\nu_o}\right)^{1/2}\dd\nu_o.
\end{equation}
Cauchy--Schwarz inequality and symmetry of the last integral in the two
probability measures give
\begin{equation}\label{eq:positive}
    0<\ground(\lambda)\leq1,\quad \ground(0)=1,\quad
    \ground(\lambda^*)=\ground(\lambda).
\end{equation}
From \eqref{eq:boundary} we may conclude that
\begin{equation}
    \label{eq:bound-by-ground}
    |\varphi_{\iu \theta}(\lambda)| \leq \varphi_0(\lambda),
    \qquad \theta \in \mfrak a, \lambda \in P^+.
\end{equation}
For $z\in\mathfrak a_\C$ and $\lambda \in P^+$,
\[
    A_\lambda\varphi_z=P_\lambda(z)\varphi_z.
\]
In particular,
\begin{equation}\label{eq:eigen}
    A_\lambda\varphi_0=\ground(\lambda)\varphi_0
    \qquad \lambda \in P^+.
\end{equation}
A proof can be found in \cite[Theorem 3.22(iii)]{ParkinsonSHA}.
We will also use the following identity
\begin{equation}\label{eq:sphere-plane-wave}
    \sum_{y\in V_\lambda(x)}
    \chi_0(h(x,y;\omega))^{1/2}\e^{\langle z,h(x,y;\omega)\rangle}
    =N_\lambda P_\lambda(z),
    \qquad x \in V_g, \lambda \in P^+, \omega \in \Omega,z \in \mfrak a_\C,
\end{equation}
which is \cite[Corollary 3.23 and Theorem 6.2]{ParkinsonSHA};
see also \cite[Proposition 4.1]{PapageorgiouTrojan}.
These statements do not require $\varphi_z\in\ell^2(V_g)$.

\subsection{Integrability of the ground spherical function}\label{sec:integrability}
The estimate underlying the integrability argument is
\begin{equation}\label{eq:ground-estimate}
    \ground(\lambda)\approx\chi_0(\lambda)^{-1/2}
    \prod_{\alpha\in\Phi^{++}}(1+\langle\alpha,\lambda\rangle).
\end{equation}
For the building statement matching our setting see
\cite[Proposition 1.2]{PapageorgiouTrojan}; for type $\widetilde A_r$,
see \cite[Proposition 2.1]{AnkerSchapiraTrojan}.
More detailed building asymptotics
are given in \cite[Section 7.2]{RemyTrojan}.
Here $a(\lambda)\approx b(\lambda)$ for functions
$a, b\colon P^+ \to (0, \infty)$ means $a(\lambda) \lesssim b(\lambda) 
\lesssim a(\lambda)$ and $a(\lambda) \lesssim b(\lambda)$ means
existence of positive constant $C$ independent of $\lambda$
such that $a(\lambda) \leq C b(\lambda)$.

For completeness, we provide the proof of the following
\begin{proposition}\label{prop:integrability}
    The ground
    spherical function $\varphi_0(x)$ belongs to
    $\ell^t(V_g)$ for every $2<t\leq\infty$, and $\|\varphi_0\|_\infty=1$.
\end{proposition}
\begin{proof}
    Set $d=|\Phi^{++}|$ and $q_{\min}=\min_{0\leq i\leq r}q_i>1$.
    The upper bound in \eqref{eq:ground-estimate} implies
    \begin{equation}\label{eq:decay}
        \ground(\lambda)\lesssim (1+\operatorname{ht}(\lambda))^d
        \chi_0(\lambda)^{-1/2}.
    \end{equation}
    For the second inequality, each nonzero translation $t_{\lambda_i}$
    has positive length, hence
    $\chi_0(\lambda_i)=q_{t_{\lambda_i}}\geq q_{\min}$; now use
    \eqref{eq:sphere} to obtain
    \begin{equation}
        \chi_0(\lambda)\geq q_{\min}^{\operatorname{ht}(\lambda)}.
    \end{equation}
    By first equation of \eqref{eq:sphere}, $N_\lambda\lesssim \chi_0(\lambda)$.
    For $2<t<\infty$ it follows that
    \begin{align*}
        \|\varphi_0\|_t^t
 &=\sum_\lambda N_\lambda\ground(\lambda)^t\\
 &\lesssim \sum_\lambda(1+\operatorname{ht}(\lambda))^{dt}
 \chi_0(\lambda)^{1-t/2}\\
 &\leq \sum_{j=0}^\infty\binom{j+r-1}{r-1}
 (1+j)^{dt}q_{\min}^{-(t/2-1)j}<\infty.
    \end{align*}
    The last count is the number of $r$-tuples of nonnegative integers
    with sum $j$; the exponential factor dominates the polynomial.
    Formulas~\eqref{eq:positive} handle $t=\infty$.
\end{proof}

\subsection{Spectral facts and the Helgason transform}\label{sec:spectral-prelim}
For a finitely supported function $u$ on $V_g$, use the definition of
Helgason transform of \cite[Section 4.1]{PapageorgiouTrojan}:
\begin{equation}\label{eq:Helgason}
    \mathcal H(u)(z,\omega)=\sum_x u(x)
    \chi_0(h(o,x;\omega))^{1/2}\e^{-\langle z,h(o,x;\omega)\rangle},
    \qquad z \in \mathfrak a_\C,\ \omega \in \Omega.
\end{equation}
For radial $k$, \eqref{eq:sphere-plane-wave} shows that Helgason transform
of $k$ is independent of $\omega$ and
\begin{equation}\label{eq:radial-transform}
    \mathcal H(k)(z)=\sum_\lambda N_\lambda\tilde k(\lambda)P_\lambda(-z)
    = \sum_x k(x)\varphi_{-z}(x).
\end{equation}
In particular,
\begin{equation}\label{eq:H-bound}
    \mathcal H(|k|)(0)=\sum_\lambda N_\lambda|\tilde k(\lambda)|\ground(\lambda)
    =\sum_{x}|k(x)|\varphi_0(x).
\end{equation}
For an arbitrary radial $k$ satisfying $\mathcal H(|k|)(0)<\infty$,
Tonelli's theorem and \eqref{eq:sphere-plane-wave} show that the series
\eqref{eq:Helgason} is absolutely convergent at $z=0$, for every
$\omega$, and equals \eqref{eq:radial-transform}.
This includes the case when $k \not\in \ell^1$.
The finite-support and $\ell^1$ versions appear in
\cite[Proposition 4.1 and Corollary 4.2]{PapageorgiouTrojan}.

Let $\mathscr A_0=\operatorname{span}\{A_\lambda:\lambda\in P^+\}$;
it is a commutative unital algebra, whose characters can be expressed
by Macdonald spherical functions: every multiplicative functional on $\mathscr A_0$
is given by evaluation at some $z \in \mathfrak a_\C$
\begin{equation}
    h_z(A_\lambda) = P_\lambda(z), \qquad \lambda \in P^+.
\end{equation}
Moreover, $h_z = h_{z'}$ if and only if $W z + 2 \uppi \iu L = W z' + 2 \uppi \iu L$
as shown in \cite[Theorem 3.3.12(ii)]{Macdonald1971}.
Let $\mathscr A_2$ be the closure of $\mathscr A_0$ in the space of bounded
linear operators $\ell^2(V_g) \to \ell^2(V_g)$. This is a commutative
unital $C^*$-algebra \cite[Section 4]{ParkinsonSHA}. Let us define
Plancherel measure $\pi$ and characterize the spectrum of $\mathscr A_2$.
In the standard case, $\tau_\alpha \geq 1$ for all $\alpha \in \Phi$.
In this case set
\begin{equation}
    U_0 = \{ \theta \in \mfrak a \colon (\forall \alpha \in \Phi)
        \langle \theta, \alpha^\lor\rangle \leq \uppi \}.
\end{equation}
Then for all $\theta \in U_0$ the multiplicative functional $h_{\iu \theta}$
extends to $\scr A_2$ in a continuous way. Set $\scr D = U_0$, and for
any borel subset $B \subseteq \scr D$,
\begin{equation}
    \pi(B) = \left( \frac{1}{2\uppi} \right)^r \frac{W(q^{-1})}{|W|}
    \int_B \frac{\dd \theta}{|\mathbf c(\iu \theta)|^2}.
\end{equation}
In the exceptional case suppose that $\tau_\alpha < 1$ for some $\alpha \in \Phi$.
This is only possible when $\Phi$ is the $BC_r$ root system and
$q_r < q_0$. Set $b = \sqrt{q_r / q_0}$, and for $j=1, \ldots, r$ let 
\begin{equation}
    U_j = \{z \in \mfrak a_\C \colon (\forall k \neq j)z_k \in [-\uppi, \uppi]
    \land z_j = \uppi - \iu \log b\},
\end{equation}
and $U_0 = [-\uppi, \uppi]^r$. For $j = 1, \ldots, r$ and $\theta \in U_j$ define
\begin{equation}
    \phi_j(\iu\theta) = \mathbf c(-\iu\theta)
    \lim_{t\to 0}\frac{\mathbf c(\iu \theta + t e_j)}{1 - \e^t},
\end{equation}
where $(e_j)_j$ is the standard basis.
Then for each $\theta \in \scr D = U_0 \sqcup U_1$, the multiplicative functional
$h_{\iu\theta}$ extends to $\scr A_2$ in a continuous way. For a Borel
subset $B \subseteq \scr D$ put
\begin{equation}
    \pi(B) = \left( \frac 1 {2\uppi} \right) ^r \frac{W(q^{-1})}{|W|}
    \int_{B \cap U_0} \frac{\dd \theta}{|\mathbf c(\iu \theta)|^2}
    + \left( \frac 1 {2\uppi} \right) ^{r-1} \frac{W(q^{-1})}{|W'|}
    \int_{B \cap U_1} \frac{\dd \theta}{\phi_1(\iu \theta)},
\end{equation}
where $W'$ is the Coxeter group $C_{r-1}$ and the measure $\dd \theta$
on $U_j$ is $\dd \theta = \prod_{k \neq j} \dd \theta_k$ with $k$
in range $1 \ldots r$ and $j = 0, 1, \ldots r$.
In either standard or exceptional case,
for all $A \in \scr A_0, x\in V_g$ and $y \in V_\lambda(x)$,
\begin{equation}
    (A\delta_x)(y) = \int_{\scr D}h_{\iu \theta}(A)
    \overline{P_\lambda(\iu \theta)} \pi(\dd \theta).
\end{equation}
See \cite[Theorem 5.2 and Corollary 5.5]{ParkinsonSHA} for the standard case,
and \cite[Theorem 5.7 and Corollary 5.8]{ParkinsonSHA} for the exceptional case.

When $f, g \colon V_g \to \C$ are finitely supported, then
\begin{equation}\label{th:4.3}
    \sum_{y \in V_g} f(y)\overline{g(y)} =
    \int_{\scr D}\int_{\Omega} \cal H(f)(\iu z, \omega)
        \cal H(\overline g)(-\iu z, \omega)\nu_o(\dd \omega)\pi(\dd z).
\end{equation}
For proof see \cite[Theorem 4.3]{PapageorgiouTrojan}.

Let $k \in \ell^1(V_g)$ be radial function. Then for each $f \in \ell^2(V_g)$,
we have
\begin{equation}\label{prop:4.4}
    \cal H(f \times k)(\iu z, \omega) = \cal H (k)(\iu z) \cal H (f)(\iu z, \omega)
    \qquad z \in \scr D, \omega \in \Omega.
\end{equation}
For proof see \cite[Proposition 4.4]{PapageorgiouTrojan}.

The bound \eqref{eq:bound-by-ground} can be generalized to
\begin{equation}
    \label{eq:bound_2}
    |\varphi_{\iu z}(x)| \leq \varphi_0(x),
    \qquad x \in V_g, z \in \scr D.
\end{equation}
For proof see \cite[Theorem B.3.5]{ParkinsonPhD}.

\subsection{Interpolation tools}\label{sec:interpolation-prelim}
We state versions of complex interpolation theorems sufficient
for our use in Section \ref{sec:kunze}.
On a countable set with counting measure, a simple function with support
of finite measure is precisely a finitely supported function.

\begin{theorem}[bounded Stein interpolation]\label{th:stein}
    Let $X,Y$ be countable sets with counting measure. For
    $0\leq\Re z\leq1$, let $S_z$ be a linear map from finitely supported
    functions on $X$ to functions on $Y$. Suppose that for every such
    $u$ and every finitely supported $v$ on $Y$, the scalar function
    \[
        z\mapsto\sum_{y\in Y}v(y)S_zu(y), \qquad
        0 \leq \Re z \leq 1,
    \]
    is bounded and continuous on the closed strip and analytic in its
    interior. Suppose $1\leq p_0,p_1,q_0,q_1\leq\infty$ and
    \[
        \|S_{\iu t}u\|_{q_0}\leq C_0\|u\|_{p_0},\qquad
        \|S_{1+\iu t}u\|_{q_1}\leq C_1\|u\|_{p_1}\qquad t\in\R,
    \]
    where $C_0,C_1$ are positive constants independent of $t$ and $u$.
    For $0<\theta<1$, set
    \[
        \frac1{p_\theta}=\frac{1-\theta}{p_0}+\frac\theta{p_1},\qquad
        \frac1{q_\theta}=\frac{1-\theta}{q_0}+\frac\theta{q_1},
    \]
    with the usual convention of $1/\infty = 0$.
    Then
    $\|S_\theta u\|_{q_\theta}
    \leq C_0^{1-\theta}C_1^\theta\|u\|_{p_\theta}$
    for finitely supported $u$. If $p_\theta<\infty$, the operator
    extends by density to $\ell^{p_\theta}(X)$.
\end{theorem} \noindent
This is the bounded-growth case of \cite{SteinInterpolation}.
The stronger boundedness hypothesis above suffices for all the finite
scalar sums used in this paper, so no general
admissible-growth formula is needed.

\begin{theorem}[Hadamard three-lines theorem]\label{th:3lines}
    Let $G$ be analytic on $0<\Re z<1$ and bounded and continuous on
    $0\leq\Re z\leq1$. If $|G(\iu t)|\leq M_0$ and
    $|G(1+\iu t)|\leq M_1$ for every $t\in\R$, where $M_0,M_1>0$,
    then
    \[
        |G(\theta+\iu t)|\leq M_0^{1-\theta}M_1^\theta
        \qquad(0\leq\theta\leq1,\ t\in\R).
    \]
\end{theorem} \noindent
For a proof see for example \cite[Lemma 1.1.2]{BerghLofstrom}.

\section{Radial Kunze--Stein for buildings}\label{sec:kunze}
Our main result the following
\begin{theorem}[Weighted radial estimate and Kunze--Stein]\label{thm:KS}
    Let $\mathscr X$ be a thick, locally finite, regular irreducible affine
    building. For a radial complex function
    $k(x)=\tilde k(\sigma(o,x))$ with
    $$\mathcal H(|k|)(0)=\sum_x|k(x)|\varphi_0(x)<\infty,$$ the operator
    \[
        T_ku(x)=\sum_y u(y)\tilde k(\sigma(x,y)),
    \]
    initially defined on finite-support functions, extends to $\ell^2(V_g)$
    and satisfies
    \begin{equation}\label{eq:weighted}
        \|T_ku\|_2\leq\mathcal H(|k|)(0)\|u\|_2.
    \end{equation}
    For $1\leq p<2$, radial $k\in\ell^p(V_g)$ and $u\in\ell^2(V_g)$,
    \begin{equation}\label{eq:KS}
        \|T_ku\|_2\leq\|\varphi_0\|_{p'}\|k\|_p\|u\|_2,
        \qquad 1/p+1/p'=1.
    \end{equation}
    Norms use counting measure, and $\|\varphi_0\|_{p'}$ is finite.
\end{theorem}
\subsection{Four proofs of the Kunze--Stein phenomenon}
The first proof is by far the simplest, however it requires
a nontrivial theorem (the equation \eqref{eq:5}).
\begin{proof}[first proof of Theorem \ref{thm:KS}]
    By \eqref{eq:kernel} and \eqref{eq:5},
    \begin{align}
        \norm{T_k}_{2\to 2}
        & \leq \sum_{\lambda} |\tilde k(\lambda)| N_\lambda \norm{A_\lambda}_{2\to 2} \\
        & \leq \sum_{\lambda} |\tilde k(\lambda)| N_\lambda \Phi(\lambda) \\
        & = \sum_{x} |k(x)| \varphi_0(x).
    \end{align}
    Equation \eqref{eq:KS} follows from Proposition \ref{prop:integrability}
    as well as the H\"older's inequality.
\end{proof}

\noindent In the second proof, we conjugate by the ground spherical function and interpolate.
Let $Mu(x)=\varphi_0(x)u(x)$. Positivity of $\varphi_0$ allows
defining $M^z$ for all $z \in \C$.

\begin{lemma}[Conjugated endpoint bounds]\label{lem:endpoints}
    For a radial function $k\colon V_g \to \C$ suppose that
    $$\mathcal H(|k|)(0)=\sum_x|k(x)|\varphi_0(x)<\infty.$$
    Then the following operator compositions are well-defined and satisfy
    \begin{align}
        \|M^{-1}T_kMu\|_\infty&\leq\mathcal H(|k|)(0)\|u\|_\infty
            \phantom{1_1}\qquad u\in\ell^\infty(V_g), \label{eq:endinf}\\
        \|MT_kM^{-1}u\|_1&\leq\mathcal H(|k|)(0)\|u\|_1
            \phantom{1_\infty}\qquad u\in\ell^1(V_g). \label{eq:end_1}
    \end{align}
\end{lemma} \noindent
Note that intermediate operators need not be bounded.
\begin{proof}
    For $u\in\ell^\infty(V_g)$, the eigenfunction identity \eqref{eq:eigen}
    and Helgason transform formula \eqref{eq:H-bound} give
    \begin{align*}
        \sum_y|\tilde k(\sigma(x,y))|\varphi_0(y)|u(y)|
 &\leq\|u\|_\infty\sum_\lambda
 N_\lambda|\tilde k(\lambda)|A_\lambda\varphi_0(x)\\
 &=\mathcal H(|k|)(0)\|u\|_\infty\varphi_0(x)<\infty.
    \end{align*}
    Divide by $\varphi_0(x)$ and take the supremum over $x$ to obtain
    the first endpoint \eqref{eq:endinf}. For $u\in\ell^1(V_g)$, Tonelli's theorem with
    \eqref{eq:opposition} and \eqref{eq:positive} give
    \begin{align*}
        \sum_x\varphi_0(x)\sum_y
        |\tilde k(\sigma(x,y))|\frac{|u(y)|}{\varphi_0(y)}
 &=\sum_y\frac{|u(y)|}{\varphi_0(y)}
 \sum_\lambda N_\lambda|\tilde k(\lambda)|A_{\lambda^*}\varphi_0(y)\\
 &=\mathcal H(|k|)(0)\|u\|_1<\infty.
    \end{align*}
    This proves the second endpoint estimate \eqref{eq:end_1}.
\end{proof}

\begin{proof}[second proof of Theorem \ref{thm:KS}]
    Set $S_z=M^{2z-1}T_kM^{1-2z}$. For finitely supported $u$,
    \[
        S_zu(x)=\sum_{y\in\supp u}
        \varphi_0(x)^{2z-1}\varphi_0(y)^{1-2z}
        \tilde k(\sigma(x,y))u(y).
    \]
    For each $x$, this is a finite linear combination of exponentials
    of the form $\exp(a z)$ with $a \in \R$,
    hence it is entire and bounded on
    $0\leq\Re z\leq1$. The same holds after pairing with any finitely
    supported function. Multiplication by $\varphi_0^{\iu t}$ is isometric
    on every $\ell^p$, so Lemma~\ref{lem:endpoints} gives
    \[
        \|S_{\iu t}u\|_\infty\leq\mathcal H(|k|)(0)\|u\|_\infty,
        \qquad
        \|S_{1+\iu t}u\|_1\leq\mathcal H(|k|)(0)\|u\|_1.
    \]
    Apply Theorem~\ref{th:stein} with
    $(p_0,q_0)=(\infty,\infty)$, $(p_1,q_1)=(1,1)$ and
    $\theta=1/2$. Since $S_{1/2}=T_k$, this proves
    \eqref{eq:weighted} for finite-support functions, and then by density on
    $\ell^2(V_g)$.
\end{proof}

\noindent
The third proof is the direct weighted Schur test argument.
Similarly to the second proof, it takes advantage of Lemma \ref{lem:endpoints}.
This proof does not need interpolation.
\begin{proof}[third proof of Theorem \ref{thm:KS}]
    Put
    $K(x,y)=|\tilde k(\sigma(x,y))|$. Cauchy--Schwarz inequality gives
    \[
        |T_ku(x)|^2\leq
        \left(\sum_yK(x,y)\frac{\varphi_0(y)}{\varphi_0(x)}\right)
        \left(\sum_yK(x,y)\frac{\varphi_0(x)}{\varphi_0(y)}|u(y)|^2\right).
    \]
    The first factor is $\mathcal H(|k|)(0)$ by \eqref{eq:H-bound}.
    Second factor is precisely $MT_{|k|}M^{-1} |u(x)|^2$, so
    summing over $x$ and using the identity \eqref{eq:end_1} gives
    \[
        \|T_ku\|_2^2\leq
        \bigl(\mathcal H(|k|)(0)\bigr)^2\|u\|_2^2.
    \]
    The same Schur argument for a single sphere operator appears in
    \cite[proof of Theorem 6.3]{ParkinsonSHA}.
\end{proof}

\noindent
The following proof was communicated to us by Effie Papageorgiou.
\begin{proof}[fourth proof of Theorem \ref{thm:KS}]
    Let $f$ be finitely supported. For finitely supported and radial
    $k\colon V_g \to \C$, using \eqref{th:4.3} and \eqref{prop:4.4}
    in the standard case of Section \ref{sec:spectral-prelim}.
    \begin{align}\label{eq:2}
        \|f \times k \|_2^2
        &= \int_{\scr D}\int_{\Omega}|\cal H(f \times k)(\iu z, \omega)|^2
        \nu_o(\dd \omega)\pi(\dd z) \\
        &= \int_{\scr D}\int_{\Omega}|\cal H(f)(\iu z, \omega)|^2
        |\cal H(k)(\iu z)|^2 \nu_o(\dd \omega)\pi(\dd z) \\
        &\leq \sup_{z \in \scr D}|\cal H(k)(\iu z)|^2
        \int_{\scr D}\int_{\Omega}|\cal H(f)(\iu z, \omega)|^2
        \nu_o(\dd \omega)\pi(\dd z) \\
        &= \sup_{z \in \scr D}|\cal H(k)(\iu z)|^2 \|f\|_2^2.
    \end{align}
    Since for radial $k$ we have \eqref{eq:radial-transform},
    applying \eqref{eq:bound-by-ground}
    \begin{equation}
        \label{eq:1}
        |\cal H(k)(\iu z)| \leq \sum_x |k(x)||\varphi_{\iu z}(x)|
        \leq \sum_x |k(x)| \varphi_0(x) = \cal H(|k|)(0).
    \end{equation}
    In the exceptional case we may not use \eqref{eq:bound-by-ground},
    or in fact write \eqref{eq:2}, since the spectrum is not real.
    Instead, we use \eqref{th:4.3} and \eqref{prop:4.4} and write
    \begin{align}
        \|f \times k \|_2^2
        &= \int_{\scr D}\int_{\Omega}
        \cal H(f \times k)(\iu z, \omega)
        \cal H(\overline{f \times k})(-\iu z, \omega)
        \nu_o(\dd \omega)\pi(\dd z) \\
        &= \int_{\scr D} \cal H(k)(\iu z) \cal H(\bar k)(-\iu z) \int_{\Omega}
        \cal H(f)(\iu z, \omega) \cal H(\bar f)(-\iu z, \omega)
        \nu_o(\dd \omega)\pi(\dd z) \\
        &\leq \sup_{z \in \scr D}|\cal H(k)(\iu z) \cal H(\bar k)(-\iu z)|
        \int_{\scr D} \left|\int_{\Omega}
        \cal H(f)(\iu z, \omega) \cal H(\bar f)(-\iu z, \omega)
        \nu_o(\dd \omega)\right|\pi(\dd z) \label{eq:4}\\
        &\leq \sup_{z \in \scr D}|\cal H(k)(\iu z)|^2 \|f\|_2^2.
    \end{align}
    The last inequality needs justification, why the absolute value
    can be omitted.
    For finitely supported $f$ and $z \in \scr D$, define
    $$Q_z(f) = \int_{\Omega}\cal H(f)(\iu z, \omega)
    \cal H(\bar f)(-\iu z, \omega)\nu_o(\dd \omega).$$
    We show that for almost all $z \in \scr D$, $Q_z(f) \geq 0$.
    By \eqref{eq:Helgason}, \eqref{eq:RN}, \eqref{eq:boundary} and
    the cocycle identity we have
    \begin{align}
        Q_z(f) &= \int_\Omega
        \left(\sum_x f(x) \chi_0(h(o, x; \omega))^{1/2}
            \e^{-\left<\iu z, h(o, x; \omega) \right>}\right) \\
        &\phantom{=\int}\cdot
        \left(\sum_y \overline{f(y)} \chi_0(h(o, y; \omega))^{1/2}
            \e^{\left<\iu z,  h(o, y; \omega) \right>}\right)
        \nu_o(\dd \omega) \\
        &= \sum_{x, y} f(x) \overline{f(y)} \int_\Omega
        \left( \frac{\dd \nu_x}{\dd \nu_o}(\omega) \right) ^{1/2}
        \left( \frac{\dd \nu_y}{\dd \nu_o}(\omega) \right) ^{1/2} 
        \e^{\left <\iu z, h(x, o; \omega) + h(o, y; \omega) \right>}\nu_o(\dd \omega)\\
        &= \sum_{x, y} f(x) \overline{f(y)} \int_{\Omega}
        \left( \frac{\dd \nu_y}{\dd \nu_x} (\omega)\right) ^{1/2}
        \e^{\left <\iu z, h(x, y; \omega)\right>}\nu_x(\dd \omega)\\
        &= \sum_{x, y}f(x)\overline{f(y)}P_{\sigma(x, y)}(\iu z).
        \label{eq:3}
    \end{align}
    The right-hand side of \eqref{eq:3} is nonnegative because
    $(x, y) \mapsto P_{\sigma(x, y)}(\iu z)$ is positive definite for almost
    all $z \in \scr D$. For a very general proof of this fact see
    \cite[Corrolary 4.5]{Voit2016}. This allows us to omit absolute value
    in \eqref{eq:4} and again use \eqref{th:4.3}.

    Using \eqref{eq:bound_2} in place of \eqref{eq:bound-by-ground} in
    \eqref{eq:1} proves it in the special case as well.

    Starting from general radial $k$ satisfying
    $\cal H(|k|)(0) < \infty$, denote $k_R = k \ind{B(o, R)}$
    where $B(o, R)$ is the ball of good vertices
    whose Euclidean distance to the vertex $o$ is less than $R > 0$.
    Notice that $k_R$ stay radial, and that they are finitely
    supported. Operators $T_{k_R}$ extend
    to $\ell^2(V_g)$ by density, with
    $\norm{T_{k_R}}_{2\to 2} \leq \cal H(|k_R|)(0)$.
    Hence, for $R' > R > 0$, we have
    \begin{equation}
        \norm{T_{k_R} - T_{k_{R'}}}_{2 \to 2} \leq
        \cal H(|k_R - k_{R'}|)(0) \xrightarrow{R \to \infty} 0,
    \end{equation}
    as it is bounded by the tail of convergent series
    $\cal H(|k|)(0)$. Let $T$ be the operator norm limit
    of $T_{k_R}$.
    For any $f \in \ell^2(V_g)$ and $x \in V_g$, the pointwise limit
    \begin{equation}
        T_{|k_n|} |f|(x) \uparrow \sum_y |\tilde k(\sigma(x, y))| |f(y)|
    \end{equation}
    exist. Since absolutely convergent sum is convergent, this identifies
    $T$ as $T_k$.
\end{proof}

\subsection{The mixed-exponent estimate}\label{sec:mixed}
We interpolate the weighted $\ell^2$ estimate with the elementary
$\ell^{w'}$ to $\ell^\infty$ estimate. The method is the weighted
interpolation used in \cite[Lemma 5.1]{AnkerPierfeliceVallarino};
we give the details for buildings.

\begin{theorem}[Mixed-exponent estimate]\label{thm:mixed}
    Let $\mathscr X$ be a thick, locally finite, regular irreducible affine
    building. Let $k(x)=\tilde k(\sigma(o,x))$ be a radial complex function.
    For $2\leq a\leq b<\infty$, define
    \[
        \frac1s=\frac1a+\frac1b,\qquad \theta=\frac2b,\qquad
        \mu=s\theta=\frac{2a}{a+b},\qquad \frac1a+\frac1{a'}=1.
    \]
    If $\varphi_0^\theta k\in\ell^s(V_g)$, the operator
    $T_ku(x)=\sum_yu(y)\tilde k(\sigma(x,y))$ maps
    $\ell^{a'}(V_g)$ boundedly to $\ell^b(V_g)$, with
    \begin{equation}\label{eq:mixed}
        \|T_ku\|_b\leq\|u\|_{a'}
        \left(\sum_{\lambda\in P^+}N_\lambda\ground(\lambda)^\mu
        |\tilde k(\lambda)|^s\right)^{1/s}
        =\|u\|_{a'}\|\varphi_0^\theta k\|_s.
    \end{equation}
    All norms use counting measure. No radiality assumption is made on
    $u$, and the kernel sum converges absolutely at every vertex.
\end{theorem}
\begin{proof}
    For $1\leq w\leq\infty$, H\"older in $y$ at each fixed $x$ gives
    \begin{equation}\label{eq:young}
        \|T_{k_g}u\|_\infty\leq\|k_g\|_{w}\|u\|_{w'}.
    \end{equation}
    Here $k_g(x)=g(\sigma(o,x))$ is the radial lift of the profile $g$.
    Theorem~\ref{thm:KS} gives the other endpoint
    \begin{equation}\label{eq:otherendpoint}
        \|T_{k_{\ground^{-1}g}}u\|_2\leq\|k_g\|_{1}\|u\|_2,
        \qquad (\ground^{-1}g)(\lambda)=\ground(\lambda)^{-1}g(\lambda).
    \end{equation}
    If $b=2$, then $a=2$, $s=1$, and this is \eqref{eq:weighted}.
    Otherwise $\theta=2/b\in(0,1)$; choose $w$ by
    \[
        \frac1{w}=\frac{1/a-1/b}{1-2/b}\in[0,1/2].
    \]
    Then
    \begin{equation}\label{eq:interpexponents}
        \frac1{a'}=\frac{1-\theta}{w'}+\frac\theta2,\qquad
        \frac1s=\frac{1-\theta}{w}+\theta,\qquad
        \frac1{b'}=1-\theta+\frac\theta2.
    \end{equation}
    Take finite-support $u,g,v$, with $g$ a profile on $P^+$, and normalize
    $\|u\|_{a'}=\|k_g\|_{s}=\|v\|_{b'}=1$; zero inputs are trivial.
    For $\zeta\ne0$ put $\operatorname{ph}(\zeta)=\zeta/|\zeta|$.
    On the respective supports define
    \begin{align*}
        u_z(x)&=\operatorname{ph}(u(x))
        |u(x)|^{a'((1-z)/w'+z/2)},\\
        g_z(\lambda)&=\operatorname{ph}(g(\lambda))
        |g(\lambda)|^{s((1-z)/w+z)},\\
        v_z(x)&=\operatorname{ph}(v(x))|v(x)|^{b'(1-z+z/2)},
    \end{align*}
    and set them to zero off their supports.
    Consider the bilinear pairing
    \[
        G(z)=\sum_xv_z(x)T_{k_{\ground^{-z}g_z}}u_z(x).
    \]
    It is a linear combination of exponentials with of the form
    $\exp(az), a \in \R$ so it is entire and bounded on the closed strip.
    For $z=\iu t$, the endpoint norms are
    $\|u_z\|_{w'}=\|k_{g_z}\|_{w}=\|v_z\|_1=1$.
    This includes $w=\infty$, since then $|g_{\iu t}|=1$ on its support.
    As $|\ground(\lambda)^{-\iu t}|=1$, \eqref{eq:young} gives
    $|G(\iu t)|\leq1$. On $\Re z=1$, the endpoint norms are
    $\|u_z\|_2=\|k_{g_z}\|_{1}=\|v_z\|_2=1$, so
    \eqref{eq:otherendpoint} gives $|G(1+\iu t)|\leq1$.
    Theorem~\ref{th:3lines} implies $|G(\theta)|\leq1$.
    By \eqref{eq:interpexponents}, all three exponents at $z=\theta$
    equal one, hence $u_\theta=u$, $g_\theta=g$ and $v_\theta=v$.
    Taking the supremum over finite-support $v$ with $\|v\|_{b'}\leq1$
    and restoring the normalization gives
    \[
        \|T_{k_{\ground^{-\theta}g}}u\|_b
        \leq\|u\|_{a'}\|k_g\|_{s}.
    \]
    Substitute $g(\lambda)=\ground(\lambda)^\theta\tilde k(\lambda)$
    to obtain \eqref{eq:mixed}.
    Since $s,a'<\infty$, truncation and density give the bounded extension.
\end{proof}

\bibliographystyle{abbrv}
\bibliography{affine_buildings_article}

@article{Voit2016,
    title={Generalized commutative association schemes, hypergroups, and positive product formulas},
    author={Voit, Michael}, year={2016},
    doi={10.31390/cosa.10.4.11},
    journal={Communications on Stochastic Analysis},
    pages={561–-585}, volume={10}, number={4}
}

@phdthesis{ParkinsonPhD,
    title={Buildings and {H}ecke algebras},
    author={Parkinson, J.},
    year={2005},
    school={Univ. of Sydney}
}

@article{ParkinsonSHA,
    author={Parkinson, James}, title={Spherical harmonic analysis on affine buildings},
    journal={Mathematische Zeitschrift}, volume={253}, number={3},
    pages={571--606}, year={2006}, doi={10.1007/s00209-005-0924-4},
    url={https://www.maths.usyd.edu.au/u/jamesp/2.pdf}
}

@book{Macdonald1971,
    author={Macdonald, I. G.}, title={Spherical Functions on a Group of {$p$}-adic Type},
    series={Publications of the Ramanujan Institute}, number={2},
    publisher={Ramanujan Institute, University of Madras}, address={Madras}, year={1971}
}

@article{KunzeStein,
    author={Kunze, R. A. and Stein, E. M.},
    title={Uniformly bounded representations and harmonic analysis of the {$2\times2$} real unimodular group},
    journal={American Journal of Mathematics}, volume={82}, number={1}, pages={1--62}, year={1960}
}

@article{VecaTree,
    author={Veca, Alessandro}, title={The {Kunze--Stein} phenomenon on the isometry group of a tree},
    journal={Bulletin of the Australian Mathematical Society}, volume={65}, number={1},
    pages={153--174}, year={2002}, doi={10.1017/S0004972700020177}
}

@article{AnkerPierfeliceVallarino,
    author={Anker, Jean-Philippe and Pierfelice, Vittoria and Vallarino, Maria},
    title={The wave equation on hyperbolic spaces}, journal={Journal of Differential Equations},
    volume={252}, number={10}, pages={5613--5661}, year={2012}, doi={10.1016/j.jde.2012.01.031}
}

@article{PapageorgiouTrojan,
    author={Papageorgiou, Effie and Trojan, Bartosz},
    title={Mass functions and asymptotic behavior of caloric functions on affine buildings},
    journal={Mathematische Annalen}, volume={395}, year={2026}, note={Article 30},
    doi={10.1007/s00208-026-03453-1}
}

@article{SteinInterpolation,
    author={Stein, Elias M.}, title={Interpolation of linear operators},
    journal={Transactions of the American Mathematical Society}, volume={83}, number={2},
    pages={482--492}, year={1956}, doi={10.1090/S0002-9947-1956-0082586-0}
}

@article{Herz1970,
    author={Herz, Carl}, title={Sur le ph\'enom\`ene de {Kunze--Stein}},
    journal={Comptes Rendus de l'Acad\'emie des Sciences de Paris, S\'erie A},
    volume={271}, pages={491--493}, year={1970}
}

@article{Cowling1978,
    author={Cowling, Michael}, title={The {Kunze--Stein} phenomenon},
    journal={Annals of Mathematics}, series={2}, volume={107}, number={2},
    pages={209--234}, year={1978}
}

@incollection{CowlingHerz1997,
    author={Cowling, Michael},
    title={{Herz}'s ``principe de majoration'' and the {Kunze--Stein} phenomenon},
    booktitle={Harmonic Analysis and Number Theory (Montreal, PQ, 1996)},
    series={CMS Conference Proceedings}, volume={21}, pages={73--88},
    publisher={American Mathematical Society}, address={Providence, RI}, year={1997}
}

@article{CowlingGiuliniMeda1993,
    author={Cowling, Michael and Giulini, Saverio and Meda, Stefano},
    title={{$L^p$--$L^q$} estimates for functions of the {Laplace--Beltrami} operator on noncompact symmetric spaces. {I}},
    journal={Duke Mathematical Journal}, volume={72}, number={1},
    pages={109--150}, year={1993}, doi={10.1215/S0012-7094-93-07206-7}
}

@article{CowlingMedaSetti2000,
    author={Cowling, Michael and Meda, Stefano and Setti, Alberto G.},
    title={Estimates for functions of the {Laplace} operator on homogeneous trees},
    journal={Transactions of the American Mathematical Society},
    volume={352}, number={9}, pages={4271--4293}, year={2000}
}

@article{ParkinsonHecke,
    author={Parkinson, James}, title={Buildings and {Hecke} algebras},
    journal={Journal of Algebra}, volume={297}, number={1}, pages={1--49},
    year={2006}, doi={10.1016/j.jalgebra.2005.08.036},
    url={https://www.maths.usyd.edu.au/u/jamesp/1.pdf}
}

@article{AnkerSchapiraTrojan,
    author={Anker, Jean-Philippe and Schapira, Bruno and Trojan, Bartosz},
    title={Sharp estimates for distinguished random walks on affine buildings of type {$\widetilde A_r$}},
    journal={Indagationes Mathematicae}, volume={36}, number={2},
    pages={383--412}, year={2025}, eprint={2312.01781}, archivePrefix={arXiv}
}

@article{RemyTrojan,
    author={R\'emy, Bertrand and Trojan, Bartosz},
    title={Martin compactifications of affine buildings},
    journal={Journal of the European Mathematical Society},
    volume={28}, number={5}, pages={1979--2066}, year={2026},
    doi={10.4171/JEMS/1771}
}

@article{Cartwright2001,
    author={Cartwright, Donald I.},
    title={Spherical harmonic analysis on buildings of type {$\widetilde A_n$}},
    journal={Monatshefte f\"ur Mathematik}, volume={133}, number={2},
    pages={93--109}, year={2001}
}

@book{BerghLofstrom,
    author={Bergh, J\"oran and L\"ofstr\"om, J\"orgen},
    title={Interpolation Spaces: An Introduction},
    series={Grundlehren der mathematischen Wissenschaften}, volume={223},
    publisher={Springer-Verlag}, address={Berlin--Heidelberg}, year={1976},
    doi={10.1007/978-3-642-66451-9}
}
\end{document}